\documentclass[aos,preprint,dvipsnames]{imsart}

\RequirePackage[round, authoryear]{natbib}

\RequirePackage[OT1]{fontenc}
\RequirePackage{graphicx}

\RequirePackage[cmex10]{amsmath}
\RequirePackage[round, authoryear]{natbib}
\usepackage{graphicx}
\usepackage{epstopdf}
\usepackage{amsfonts}
\usepackage{amsmath,amssymb}		
\usepackage{bbm} 	
\usepackage{amsthm}
\usepackage{xcolor}
\usepackage{accents}
\usepackage{easy-todo}
\usepackage{mathtools}
\usepackage{tikz}
\usepackage{paralist} 

\usepackage{relsize}
\usepackage{bbm}

\startlocaldefs

\numberwithin{equation}{section}

\newcommand{\dd}{\mathrm{d}}             % derivative
\newcommand{\thetab}{\boldsymbol\theta}

\usepackage{stackengine}
\stackMath
\newcommand\tenq[2][1]{%
\def\useanchorwidth{T}%
\ifnum#1>1%
\stackunder[0pt]{\tenq[\numexpr#1-1\relax]{#2}}{\scriptscriptstyle\thicksim}%
\else%
\stackunder[1pt]{#2}{\scriptscriptstyle\thicksim}%
\fi%
}

\makeatletter
\newcommand*\rel@kern[1]{\kern#1\dimexpr\macc@kerna}
\newcommand*\widebar[1]{%
  \begingroup
  \def\mathaccent##1##2{%
    \rel@kern{0.8}%
    \overline{\rel@kern{-0.8}\macc@nucleus\rel@kern{0.2}}%
    \rel@kern{-0.2}%
  }%
  \macc@depth\@ne
  \let\math@bgroup\@empty \let\math@egroup\macc@set@skewchar
  \mathsurround\z@ \frozen@everymath{\mathgroup\macc@group\relax}%
  \macc@set@skewchar\relax
  \let\mathaccentV\macc@nested@a
  \macc@nested@a\relax111{#1}%
  \endgroup
}
\makeatother

\graphicspath {{figures/}}

\endlocaldefs

\begin{document}

\begin{frontmatter}
\title{Quantiles for multivariate distribution \\ a short survey of population concepts% and Multiple-Output Quantile Regression
}
\runtitle{Multivariate Quantiles}

% \

$\,$ \\ 

%\today

\begin{aug}
\author%[A]
{\fnms{Marc}~\snm{Hallin}%\thanksref{m1}\thanksref{m2}
\ead[label=e1]
{mhallin@ulb.ac.be}}
%\author[B]{\fnms{Gilles} \snm{Mordant,}%\thanksref{m3}
%\ead[label=e2]{gilles.mordant@yale.edu}}
% and 
%\author[B]{\fnms{Miroslav} \snm{ \v{S}iman}%\thanksref{m3}
%\ead[label=e2]{MirekSiman@seznam.cz}}
\
% Affiliations: centered, wrapped, fixed width
\address%[A]
{Universit\' e libre de Bruxelles, Brussels, Belgium \\ and \\ Institute of Information Theory and Automation, Czech Academy of Sciences, Prague, Czech Republic\\ \printead[presep={\ }]{e1} }
%{\thanksmark{m1}\parbox[t]{0.6\textwidth}{\centering \it 
%Université libre de Bruxelles, Belgium  \\  and \\ 
%Institute of Information Theory and Automation,
%Czech Academy of Sciences, Prague, Czech Republic\smallskip }}

%\address[B]{Department of Applied and Computational Mathematics,
%Yale University\printead[presep={,\ }]{e2}} 

%
%\address[B]{Institute of Information Theory and Automation, Czech Academy of Sciences, Prague, Czech Republic\printead[presep={\ }]{e2}} 

\runauthor{Marc Hallin}% and \v{S}iman}
\end{aug}

\begin{abstract} 
{\bf Abstract.\ } Quantiles are among the most fundamental concepts in Probability and Statistics,   from descriptive to
inferential.   However, quantile functions are  well-defined and well-understood  in the context of  one-dimensional proba\-bility distributions, where their definition as the inverse of distribution functions---a definition which is intimately related to the canonical ordering of the real line $\mathbb R$. Starting with dimension $d=2$,  such a canonical ordering is no longer avai\-lable in~${\mathbb R}^d$; nor is it  available in nonlinear manifolds such as hyperspheres, tori, or polyspheres (vectors of directional variables).  This results in the absence of an obvious and widely accepted definition of quantiles.  The need to extend the concept of quantile beyond the classical univariate context nevertheless has sparked a large body of literature,  giving rise to a variety of more or less satisfactory quantile concepts; the recent years have been particularly active in this respect.   This short review is an attempt to complement and update Serfling's 25-years old survey \citep{Serfling08} by summarizing  and unifying some of these concepts in the light of recent contributions. 
\end{abstract}

%\printaddresses

%\begin{keyword}[class=MSC]
%\kwd[Primary ]{62G10}\kwd{62G20}
%\kwd[; secondary ]{62P20}\kwd{62M10}
%\end{keyword}

%\begin{keyword}
%\kwd{semiparametric efficiency}
%\kwd{maximal ancillarity}
%\kwd{elimination of nuisances}
%\kwd{Gaussian shift}
%\kwd{Brownian drift}
%\kwd{measure transportation}
%\kwd{distribution-freeness}
%\kwd{center-outward ranks and signs}
%\end{keyword}

\end{frontmatter}

%%%%%%%%%%%%%%%%%%%%%%%%%%%%%%%%%%%%%%%%%%%%%%%%%%%%%%%%%%%%%%%%%%%%%%%%%%%%%%

%\bigskip

\section{Introduction} Quantiles are among the most fundamental concepts in Probability and Statistics,   from descriptive to
inferential.  Via quantile regression \citep{KoenkerBassett78,HallinSimanHandbook} and quantile autoregression \citep{koenker2006quantile}, and, more recently,   \citep{Barrio25,  Gonzetal26},  they moreover provide the most informative solution to regression and autoregression problems---the analysis of the impact on some variable of interest of a set of covariates or past values. Quantiles also constitute one of the most fundamental notions in survival analysis and  (value at risk, 
expected shortfall, ... )  risk analysis.% \citep{Beirlant20}. 

Now, quantile functions are a well-defined and well-understood concept in the context of univariate real-valued random variables and one-dimensional proba\-bility distributions, where their definition as the inverse of distribution functions, however,  is intimately related to the canonical ordering of the real line $\mathbb R$. Starting with dimension $d=2$,  such a canonical ordering is no longer avai\-lable in~${\mathbb R}^d$; nor is it  available in nonlinear manifolds such as the~$d$-sphere~${\mathcal S}_d$ (directional variables), the~$d$-torus ${\mathcal T}_d$ ($d$-tuples of circular variables), or polyspheres (vectors of directional variables). This results in the absence of an obvious and widely accepted definition of quantiles.  

However, the need to extend the concept of quantile beyond the classical univariate context has been strongly felt for more than half a century, sparking a large body of literature  and giving rise to a variety of quantile concepts and methods that are more or less satisfactory. As a result, as mentioned by Serfling in his influential 2002 survey  \citep{Serfling08}, the term {\it ``quantile''} has unfortunately become somewhat loosely used. 

  This short review is an attempt to complement and update Serfling's 25-years old one by summarizing  and unifying some of these concepts in the light of recent contributions; the recent years, indeed,  have been particularly active in the area.
%, emphasizing their respective advantages and disadvantages.
 Due to strict constraints on the number of pages, we limit ourselves to a discussion of population concepts, leaving aside their empirical counterparts and the related inference problems. For the same reason, our list of references is limited, hence  involves  some unavoidable  personal choices. We apologize in advance for any relevant publications that may be omitted.  
 
Throughout, for the sake of simplicity in exposition, we tacitly assume that all %variables and
 distributions considered here belong to the class $\mathcal P$ of Lebesgue-absolutely continuous\footnote{For Riemannian manifolds, the Lebesgue measure is defined as the measure associated with the canonical Riemannian metric.} distributions with densities that are bounded  and bounded away from zero on compacts. Mild as it is, this restriction is unnecessarily strict\footnote{See, e.g., \cite{delB20, delB24}.} for most of the results below; however, it has the advantage of being sufficient for all of them, throughout the survey. 
%%%%%%%%%%%%%%%%%%%%%%%%%%%%%%%%%%%%%%%%%%%%%%%%%%%%%%%%%%%%%%%%%%%%%%%%%%%%%%
\section{Univariate distributions}% and quantile functions, ranks, and signs}
\subsection{The classical definition}
The quantile function of a real-valued random variable~$Y$~with Lebesgue-absolutely continuous  distribution ${\rm P}_Y\in{\mathcal P}$ and distribution function~$F_Y : {\mathbb R}\to (0,1), y\mapsto F_Y(y)\coloneqq {\rm P}_Y[(-\infty, y]] $  %$ where $F_Y: {\mathbb R}\to[0,1], \ y\mapsto F_Y(y) \coloneqq {\rm P}[Y]$
 is defined\footnote{Under the assumption made of a nonvanishing density, $F_Y$ is strictly monotone increasing on the support of~${\rm P}_Y$ and a homeomorphism, so that $F_Y^{-1}$ is well defined, continuous, and monotone increasing. This nonvanishing assumption, of course, can be waived;  $Q_Y\coloneqq F_Y^{-1}$ then still can be   defined as a {\it generalized inverse}.} as $$Q_Y: (0,1)\to{\mathbb R}, \ u\mapsto Q_Y(u)\coloneqq F_Y^{-1}(u).$$
  A coherent extension should thus provide simultaneous definitions, as the inverse of each other,  of both quantile and distribution functions.\smallskip

Several other characterizations of $F_Y$ and $Q_Y$ are possible, though:\vspace{-1mm} 
\begin{enumerate}
\item[(i)] ($F_Y$ as a monotone probability integral transformation) $Y\mapsto F_Y(Y)$   is the unique\footnote{This follows as a particular case of a famous result by \cite{McCann95}. Uniqueness actually is only a.s.~uniqueness, that is, only holds up to a set of $Y$ values contained in a set of~${\rm P}_Y$ probability zero; in the sequel, for the sake of simplicity, we omit to repeatedly mention this fact.} {\it monotone probability integral transformation}---namely, $F_Y$ is the unique monotone increasing (the derivative of a convex function) mapping such that $F_Y(Y)\sim {\rm U}_{[0,1]}$, where~${\rm U}_{[0,1]}$ denotes the uniform distribution over $[0,1]$; equivalently, using the terminology and convenient notation of measure transportation, $F_Y$ is~{\it pushing~${\rm P}_Y$ forward to} ${\rm U}_{[0,1]}$, which we denote as~$F_Y\# {\rm P}_Y =~\!{\rm U}_{[0,1]}$; 
\item[(ii)]($Q_Y$ as a monotone quantile function) similarly, $Q_Y$ is the unique$^3$ monotone increasing  (the~derivative of a convex function)  {\it pushing the uniform~${\rm U}_{[0,1]}$ forward to~${\rm P}_Y$}, i.e.,  such that~$Q_Y\# {\rm U}_{[0,1]}\! =~\!\!{\rm P}_Y$;
\item[(iii)] (a collection of weighted $L_1$ location parameter) $Q_Y=F^{-1}_Y$   is the minimizer of the (convex in $y\in{\mathbb R}$) expected check
function: for any $\tau\in[0,1]$,\vspace{-2mm}  \end{enumerate}
\begin{align}\label{rhoquant}
Q_Y(\tau) &=\arg\! \min_{y\in{\mathbb R}} {\rm E}\big[\tau\vert Y-y\vert I[Y\leq y] + (1-\tau)\vert Y-y\vert I[ Y>y]\big]
\nonumber \\
&\eqqcolon\arg\! \min_{y\in{\mathbb R}} {\rm E}\big[\rho_\tau (Y-y)
\big];
%\nonumber  \\  
%&= \arg\!\min{}_{y\in{\mathbb R}}{\rm E}\big[\vert Y-y\vert + (2\tau -1)(Y-y)
%\big];
\end{align}
\begin{enumerate}\item[(iv)] (relation to copulas) $F_Y$ is the copula transform of~$Y$;
\item[(v)] (relation to depth) $\big(1-\vert 2F_Y(y ) -1\vert \big)/2 = \min (F_Y(y), 1-F_Y(y))$ is the Tukey depth of~$y$ with respect to ${\rm P}_Y$ (other depth concepts, with appropriate transformations to $F_Y$, also can be considered).
%\item[(i)]
%\item[(i)]
\end{enumerate}

All of these characterisations can be used to provide equivalent definitions of $F_Y$ and $Q_Y$ in dimension $d=1$. Extending them to dimensions $d\geq 2$, however, yields distinct concepts, with distinct properties: measure-transportation-based quantiles for (i)--(ii), spatial and Oja quantiles for (iii), depth-based quantiles for (v) ... The question thus naturally arises: which properties of the univariate concept do we expect from their multivariate extensions?  

Associated with $Q_Y$ are  {\it quantile regions ${\mathbb C}_Y(\tau)$ of order $\tau$}, $\tau\in (0,1)$. Such regions can be obtained\footnote{The classical definition  of the quantile region of order $\tau$, viz.\ ${\mathbb C}_Y(\tau) \coloneqq (-\infty,Q_Y(\tau)$, is a direct one and is not   explicitly based on a characterization of the form 
  $Q_Y({\scriptstyle{{\mathbb C}}}(\tau))$. For ${\scriptstyle{\mathbb C}}(\tau)= (-\infty,\tau]$, however, they clearly are equivalent, while the second option allows for more flexible choices of ${\scriptstyle{{\mathbb C}}}(\tau)$ (see examples below). This flexibility is unimportant in dimension one, but will be   useful in ${\mathbb R}^d$ $(d>1)$ and crucial on nonlinear manifolds.} as the images by $Q_Y$ of a collection (that does not depend on ${\rm P}_Y$) of nested (as~$\tau$ increases) subsets ${\scriptstyle{\mathbb C}}(\tau)$ of $(0,1)$ with, irrespective of ${\rm P}_Y\in{\mathcal P}$,  $(F_Y\#{\rm P}_Y)$-probability content $\tau$. The existence of such regions  for all $\tau$  requires the distribution-freeness of $F_Y(Y)$ (which here is uniform  over $[0,1]$, 
irrespective of ${\rm P}_Y\in{\mathcal P}$). 
 The classical choice for ${\scriptstyle{\mathbb C}}(\tau)$ is the interval $(0, \tau]$, yielding
 the one-sided  half-lines ${\mathbb C}_Y(\tau) = (-\infty, Q_Y(\tau)]$ as quantile regions,   with the property that, for all ${\rm P}_Y\in{\mathcal P}$ and $\tau\in(0,1)$,   
 $${\rm P}_Y\!\big[ {\mathbb C}_Y(\tau) \big]\!=\! {\rm P}_Y\!\big[(-\infty,Q_Y(\tau)]\big] \!=\! {\rm P}_Y\!\big[Y\leq Q_Y(\tau)\big] \!=\!{\rm P}_Y\!\big[F_Y(Y)\leq  \tau\big]\!=\!{\rm U}_{[0,1]} ((0,\tau])\!=\!~
 \!\tau .$$

 Other choices for ${\scriptstyle{\mathbb C}}(\tau)$ are possible, though, such as  ${\scriptstyle{\mathbb C}}(\tau)=[1-\tau , 1)$, yielding  (as in \cite{Loeve63}) quantile regions~${\mathbb C}_Y(\tau) = [Q_Y(1-\tau), \infty)$, or~${\scriptstyle{\mathbb C}}(\tau)=[(1-\tau)/2, (1+\tau)/2]$, with quantile regions the interquantile intervals
 $${\mathbb C}_{Y\pm}(\tau)\coloneqq [Q_Y((1-\tau)/2), Q_Y(((1+\tau)/2))].$$%\vspace{2mm} 

\subsection{Center-outward distribution and quantile functions}Before proposing such a list, however, let us highlight a lack of symmetry in the traditional univariate concept which, minor as it is, hinders its direct extension to dimensions~$d=2$ and higher. When the distribution~${\rm P}_Y$ is symmetric, i.e., when ${\rm P}_Y = {\rm P}_{-Y}\eqqcolon{\rm P}$, one would expect that~$Q_{-Y} =~\!Q_Y  = Q$, where $Q$ is the quantile function associated with the common distribution $\rm P$ of $Y$ and $-Y$. Instead,  $Q_{-Y} = 1-Q_Y$ since %, due to the fact that
 $F_Y(y)={\rm P}_{Y}[(-\infty,\,  y]]$ for all $y$ while~$F_{-Y}(y)={\rm P}_{Y}[[y,\, \infty)]$ and  the same asymmetry is affecting  quantile regions.  The reason for this is that~$-\infty$ and~$+\infty$, in the classical definition, do not play symmetric roles. Such asymmetry, in dimension~$d=2$, where each direction has a point at infinity which is neither $-\infty$ nor~$+\infty$, does not make sense anymore. Therefore, we introduce the so-called {\it center-outward distribution} and {\it quantile} functions~$F_{Y\pm}\coloneqq 2F_Y-1$ and~$Q_{Y\pm}\coloneqq F_{Y\pm}^{-1}$, respectively.  
 
 While $F_Y$ and $Q_Y$ are mappings to and from the unit interval $(0,1)$, the center-outward distribution and quantile functions $F_{Y\pm}$ and $Q_{Y\pm}$ are mappings  (transports) to and from~$(-1,1)\eqqcolon {\mathbb S}_1$, the open unit ball in $\mathbb R$. Clearly, $F_{Y\pm}$ and $Q_{Y\pm}$  both carry the same information on~${\rm P}_Y$ as $F_Y$ and $Q_Y$, and enjoy the following essential properties.
\begin{enumerate}
\item[(a)] Both $F_{Y\pm}$ and $Q_{Y\pm}$ characterize ${\rm P}_Y\in{\mathcal P}$;
%\item[(ii)] (monotonicity) both $F_Y$ and $Q_Y$ are monotone increasing, that is, the derivative (the gradient) of a convex function; 
\item[(b)] (monotonicity of $y\mapsto F_{Y\pm}(y)$ and distribution-freeness of $F_{Y\pm}(Y)$) $Y\mapsto F_{Y\pm}(Y)$ is a monotone proba\-bility integral transformation---namely,  (b1) $F_{Y\pm}$ is continuous and monotone increasing (the derivative of a convex function), (b2) it~pushes ${\rm P}_Y$ forward to a distribution that does not depend on ${\rm P}_Y$, and (b3) that distribution is  the uniform ${\rm U}_1$  over $(-1,\, 1)$ (the open unit ball ${\mathbb S}_1$ in $\mathbb R$), that is,   $F_{Y\pm}\# {\rm P}_Y =~\!{\rm U}_1$; 
\item[(c)] $Q_{Y\pm}$ (c1) % is a monotone proba\-bility integral transformation: $Y\mapsto F_Y(Y)$   is a {\it monotone proba\-bility integral transformation}---namely, $F_{Y\pm}$
 is monotone increasing (the derivative of a convex function) and~(c2) pushes the uniform~${\rm U}_1$ over ${\mathbb S}_1=(-1,\, 1)$  forward to ${\rm P}_Y$, that is,  $Q_{Y\pm}\#{\rm U}_1 =~ {\rm P}_Y$.
% where ${\rm U}_1$ is the uniform over $(-1,\, 1)$, the open unit ball ${\mathbb S}_1$ in $\mathbb R$.  
%\item[(f)] (the ``$Q\circ F$ property'') For any $Y\sim{\rm P}_Y\in{\mathcal P}$ and $Z\sim{\rm P}_Z\in{\mathcal P}$, with distribution and quantile functions $F_Y$, $Q_Y$, $F_Z$, and $Q_Z$, respectively, $(Q_Z\circ F_Y) \# {\rm P}_Y = {\rm P}_Z$.
\end{enumerate}
Note that, for ${\rm P}_Y\in{\mathcal P}$, (b) entails (c), and that  (b2)   plays a fundamental role in (d) the possi\-bility of defining  quantile regions ${\mathbb C}_{Y\pm}(\tau)$ of order $\tau$ with prescribed ${\rm P}_Y$~probability content~$\tau$:\footnote{This control over the probability of quantile regions, somehow, is the essence of quantiles: who would call {\it median} of ${\rm P}_Y$ a quantity $q_{1/2}$ such that ${\rm P}_Y[Y\leq q_{1/2}]= 0.5$ for  ${\rm P}_Y={\rm P}_1$ and ${\rm P}_Y[Y\leq q_{1/2}]= 0.6$ for~${\rm P}_Y=~\!{\rm P}_2$? or {\it interquartile interval} an interval with  ${\rm P}_Y$~probability content 0.4 for  ${\rm P}_Y={\rm P}_1$ and 0.6 for~${\rm P}_Y={\rm P}_2$?} without distribution-freeness, collections of subsets~${\scriptstyle{\mathbb C}}(\tau)$ with probability $\tau$ that do not depend on~${\rm P}_Y$ typically do not exist, while non-monotone mappings $Q_{Y\pm}$ would not preserve the probability contents of such regions.  
%
%\item[(c)] (control over the probability contents of quantile regions) %for all $0<\tau_1 \leq \tau_2<1$, 
%$${\rm P}_Y\big[Q([\tau_1,\,\tau_2])\big]= \tau_2 - \tau_1,$$
 %irrespective of ${\rm P}_Y$---implying, in particular, that 
% the {\it quantile region of order  $\tau$} has probability content $\tau$ for all $\tau\in(0,1)$ and ${\rm P}_Y$:
% ;
%\item[(f)] 
%\item[(d)] 
%
%Center-outward quantile regions. 

Now, the uniform ${\rm U}_1$ over ${\mathbb S}_1$ is invariant under symmetry with respect to the origin (which, along with identity, constitutes the group of orthogonal transformations of $\mathbb R$). Therefore, when choosing a collection  of nested subsets of ${\mathbb S}_1$  (the quantile regions of order $\tau$ for ${\rm U}_1$) to define the quantile regions of ${\rm P}_Y$, it is reasonable to restrict to symmetric (orthogonal-invariant) ones, and the symmetric central intervals $[-\tau, \tau] = \tau\overline{\mathbb S}_1$ (the nested closed balls with radii $\tau$, hence ${\rm U}_1$~probability contents $\tau$) are natural candidates. The collection, for~$\tau$ ranging over $(0,1)$, of the  images by $Q_{Y\pm}=F_{Y\pm}^{-1}$ of these nested balls  is the collection of nested interquantile intervals %of the form  
 ${\mathbb C}_{Y\pm}(\tau)\coloneqq [Q_Y((1-\tau)/2), Q_Y(((1+\tau)/2))].$ 
Call ${\mathbb C}_{Y\pm}(\tau)$ the {\it center-outward quantile region of order $\tau$ of} ${\rm P}_Y$ and define the  {\it center-outward quantile contour of order $\tau$} as its (two-point) boundary ${\mathcal C}_{Y\pm}(\tau)\coloneqq \{Q_Y((1\pm\tau)/2)\} = Q_{Y \pm}(\tau{\mathcal S}_1)$  (with ${\mathcal S}_1$  the unit sphere in $\mathbb R$).  

Easy computation moreover  shows that, substituting $({1\pm\tau})/{2}$ for~$\tau$ in   \eqref{rhoquant}, these boundary points~$Q_{Y\pm}(\pm\tau)=Q_Y((1\pm\tau)/2)$  of ${\mathbb C}_{Y\pm}(\tau)$ are obtained~as\footnote{In order to avoid moment assumptions, $\arg\! \min_{y\in{\mathbb R}} {\rm E}\big[ \vert Y-y\vert +\tau v (Y-y)\big]$ (involving a possibly infinite expectation) in \eqref{Chaudhu1} safely can be replaced with\vspace{1mm}  $\arg\! \min_{y\in{\mathbb R}} {\rm E}\big[ \vert Y-y\vert  +\tau v (Y-y) -\vert Y\vert- \tau v Y\big]$, where the expectation is always finite; subtracting $\vert Y\vert +\tau v Y$, indeed,  has no impact on the $\arg\! \min_{y\in{\mathbb R}}$. However, for clarity of exposition, we will stick to the simpler form, \eqref{Chaudhu1}.}
\begin{align}\label{Chaudhu1} 
Q_{Y\pm}(\tau v)&=
%Q_{Y\pm}(\pm\tau)&=
\arg\! \min_{y\in{\mathbb R}} {\rm E}\big[\rho_{(1+\tau v)/2} (Y)\big] \nonumber  \\ 
%\arg\! \min_{y\in{\mathbb R}} {\rm E}\big[\rho_{(1\pm \tau)/2} (Y)\big] \\ 
&= \arg\! \min_{y\in{\mathbb R}} {\rm E}\big[ \vert Y-y\vert +\tau v (Y-y)\big]\quad v\in\{-1,1\}={\mathcal S}_0,\  \tau\in[0,1)
%&= \arg\! \min_{y\in{\mathbb R}} {\rm E}\big[ \vert Y-y\vert  \pm\tau (Y-y)\big]\quad \tau\in[0,1).
\end{align}
and indexed by $u=\tau v\in (-1,1)={\mathbb S}_1$.
%The traditional choice, when $F_Y$ and $Q_Y$ are mappings to and from the unit interval $(0,1)$, is the image by $Q_Y$ of intervals of the form $(0,\tau]$, yielding the halflines ${\mathbb C}_Y \coloneqq (-\infty, Q_Y(\tau)]$ and the contours ${\mathcal C}_{Y}(\tau)\coloneqq \{-\infty,\, Q_Y(\tau)\}$, $\tau\in (0,1)$.\medskip 

With these definitions, we have the following further properties: 
\begin{enumerate}
\item[(d)] (control over the probability contents of quantile regions)  the ${\rm P}_Y$~probability content of~${\mathbb C}_{Y\pm}(\tau)$ (of ${\mathbb C}_Y(\tau)$), $\tau\in (0,1)$ is $\tau$, irrespective of ${\rm P}_Y$---a direct but crucial consequence of the distribution-freeness property (b2) of $F_{Y\pm}(Y)$ (of $F_{Y}(Y)$);
\item[(e)] (equivariance of quantile regions and quantile contours under shift, rescaling, and symmetry)\footnote{In a multivariate context,  this property takes the form of equivariance under shift, global rescaling, and orthogonal transformations---not the more general affine ones. Affine equivariance can be imposed via the classical transformation-retransformation technique \citep{Chakra96, Chakra01} but generally requires moment assumptions. We argue, however, that affine invariance is not necessarily a desirable property in dimension $d\geq 2$: see Section~\ref{SecAff} and the online appendix A6 of \citep{HHH} for a discussion. } for any $\tau\in (0,1)$,  any~$a\in{\mathbb R}$ and~$b\in{\mathbb R}\setminus\{0\}$, $${\mathbb C}_{(a+ bY)\pm }(\tau) = a + b\,{\mathbb C}_{Y\pm}(\tau)\quad\text{  and }\quad {\mathcal C}_{(a + bY) \pm }(\tau) = a + b\,{\mathcal C}_{Y \pm}(\tau).$$  
%\item[(f)] (the ``$Q\circ F$ transport property'') for any $Y\sim{\rm P}_Y\in{\mathcal P}$ and $Z\sim{\rm P}_Z\in{\mathcal P}$, with distribution and quantile functions $F_{Y\pm}$, $Q_{Y\pm}$, $F_{Z\pm}$, and $Q_{Z\pm}$, respectively, 
%$$(Q_{Z\pm}\circ F_{Y\pm}) \# {\rm P}_Y = {\rm P}_Z.$$ 
\end{enumerate}
%This property (f), which follows from (b)--(c), is important for applications in rank-based inference, as it guarantees the existence of a unique orbit 
% $$\{ (Q_Z\circ F_Y)\#{\rm P}_Y\vert\, Q_Z \text{ the quantile function of some ${\rm P}_Z\in{\mathcal P}$}\}={\mathcal P}$$ along which the ranks of an i.i.d. sample from   ${\rm P}\in{\mathcal P}$ are maximal invariant---which in turn implies the distribution-freeness of the latter.

Properties (a)--(c), hence also (d)--(e), are the properties we would like to see in satisfactory multivariate extensions  ${\bf F}_{{\bf Y}\pm}$ and ${\bf Q}_{{\bf Y}\pm}$ of ${F}_{{Y}}$ and ${Q}_{{Y}}$ (equivalently, of ${F}_{{Y}\pm }$ and ${Q}_{{Y}\pm }$).  

\subsection{Affine equivariance}\label{SecAff}It is often argued that any concept of multivariate quantile should enjoy {\it affine} equivariance instead of equivariance under shift, global rescaling, and orthogonal transformations. Although coinciding in dimension $d=1$, the two equivariance properties   are quite distinct (the second one being  much weaker) in dimension $d\geq 2$. \cite{Serfling08} makes it his first requirement for a sensible concept of multivariate quantile, while \cite{DPdNagy} call it ``not negotiable.''  

The argument that the analysis should not depend on the choice of marginal measurement units, at first sight, sounds quite compelling: measuring a distance in meters or kilometers, indeed,  should not have any impact on the analysis and, under the assumption of finite second-order moments, this is easily taken care of by marginal standardization. Now, what is it that should make {\it linear} marginal transformations so important? When quantifying a given phenomenon, statisticians and data analysts are often facing the problem of choosing between various measurement scales that are not affine transforms of each other. Acoustic power can be measured in  watts/m$^2$ or in decibels, which are related in a logarithmic way;   concentrations in H$^+$ ions are expressed in moles per volume or via pH, which are not linear functions of each other;  the luminosity of a star is measured as its radiated electromagnetic energy per unit time (in joules per second or watts), or as its magnitude (a logarithmic measure of  luminosity within some specific wavelength range); the strength of an earthquake is expressed as the amplitude of the seismic waves or in terms of Richter magnitude, the   value of a stock as its current stock price or its  call option price, ...  All these measurements are nonlinear monotone functions of each other, and do not commute with linear marginal transformations: standardizing them still yields distinct results. Yet, the same argument that measurement units should not matter applies: ideally, we could dream of equivariance under any marginal-order-preserving transformation. But that would imply the quantile function to be measurable with respect to the copula of the underlying distribution, which has (see Section~\ref{Sec31}) some less satisfactory consequences. 

Equivariance under shift, global rescaling, and orthogonal transformations, which is much weaker than affine-equivariance, only deals with   transformations, which  ``move around''   distribution swithout really ``modifying''~them.

\section{Multivariate distributions}\label{Sec3}% and quantile functions (population concepts)}
\subsection{Classical definitions and copula transforms}\label{Sec31} The traditional definition of the distribution function $ F_{\bf Y}$ %(also called {\it cumulative} distribution function)
 of a $d$-dimensional random vector~${{\bf Y} = (Y_1,\ldots,Y_d)}^\top$ with distribution~${\rm P}_{\bf Y}\in{\mathcal P}$ over ${\mathbb R}^d$ is %the  mapping 
$$F_{\bf Y}: {\mathbb R}^d\to (0,1),\ {\bf y}=(y_1,\ldots,y_d)\mapsto F_{\bf Y}({\bf y}) \coloneqq {\rm P}_{\bf Y}\big[Y_1\leq y_1,\ldots,Y_d\leq y_d
\big].$$ 

Denoting by $F_1,\ldots,F_d$ the marginal distribution functions (viz., $F_i\!\coloneqq\! F_{Y_i}$, $i=1,\ldots,d$) of $\bf Y$, this definition is closely related to the {\it copula transform} of $\bf Y$, namely, the mapping~${\bf F}_{{\bf Y}\text{\rm cop}}: ~\!{\mathbb R}^d \to (0,1)^d,$ ${\bf y}=(y_1,\ldots,y_d)\mapsto {\bf F}_{{\bf Y}\text{\rm cop}}({\bf y})\coloneqq  \big(F_1(y_1),\ldots,F_d(y_d)\big)$. This co\-pula transform $\bf F_{{\bf Y}\text{\rm cop}}$ is pushing ${\rm P}_{\bf Y}$ forward to a distribution ${\rm P}_{{\bf Y}\text{\rm cop}}$ over $[0,1]^d$ with uniform (over~$[0,1]$) margins---$\bf Y$'s {\it copula} which, however, is not uniform over $[0,1]^d$ since the marginals of $\bf Y$, typically, are not mutually independent. %, and depends on ${\rm P}_{\bf Y}$. 

For $d=1$,  ${\bf F}_{{\bf Y}\text{\rm cop}}$ reduces to the traditional univariate distribution function~$F_Y$, and is the natural extension resulting from $F_Y$'s characterization (iv), with range the open unit cube~$(0,1)^d$ instead of the interval $[0,1]$. As a multivariate distribution function, however,~${\bf F}_{{\bf Y}\text{\rm cop}}$ fails to satisfy the distribution-freeness property (b); as a consequence, the quantile regions  associated with ${\bf F}_{{\bf Y}\text{\rm cop}}^{-1}$ cannot 
%yield   quantile regions with prescribed probability contents $\tau$ irrespective of ${\rm P}_{\bf Y}$
satisfy the essential property~(d):  more precisely, there exists no collection of nested regions ${\scriptstyle {\mathbb C}}(\tau)$  %${\mathbb B}(\tau)\subset
in $ [0,1]^d$, $\tau\in (0,1)$ that do not depend on~${\rm P}_{\bf Y}$ and yet have probability contents $\tau$   for any ${\rm P}_{\bf Y}\in{\mathcal P}$.   Copula transforms and their inverses, therefore, do not qualify as valid multivariate extensions of univariate distribution and quantile functions. 

Turning to $F_{\bf Y}$, it is easy to see that  $F_{\bf Y}({\bf y})={\rm P}_{{\bf Y}\text{\rm cop}}\big[(0,F_1(y_1)]\times\ldots\times (0,F_d(y_d)]\big]$. Hence, $F_{\bf Y}({\bf Y})$  generally is not distribution-free.  Its inverse  ${\bf Q}_{{\bf Y}}\coloneqq {F}_{{\bf Y}}^{-1}$ (reducing to the traditional quantile function for $d=1$) is   set-valued, of the form 
$${\bf Q}_{{\bf Y}} (\tau) = F^{-1}_{\bf Y}(\tau) = \big\{
{\bf y}=(y_1,\ldots,y_d)\in{\mathbb R}^d \vert \,  {\rm P}_{{\bf Y}\text{\rm cop}}\big[(0,F_1(y_1)]\times\ldots\times (0,F_d(y_d)]
\big]	=\tau \big\}.
$$
Again, no collection   $\{{{\ensuremath{\mathsmaller{\mathbb{C}}}}}(\tau)\subset (0,1)\vert\, \tau\in (0,1)\}$ independent of ${\rm P}_{\bf Y}$ exists with the property that ${\rm P}_{\bf Y}\big[{F_{\bf Y}({\bf Y})\in{{\ensuremath{\mathsmaller{\mathbb{C}}}}}(\tau)}\big] = \tau$ irrespective of~${\rm P}_{\bf Y}$.  Therefore, no collection of quantile regions (of the form ${\bf Q}_{{\bf Y}} ({\ensuremath{\mathsmaller{\mathbb{C}}}}(\tau))$ can have ${\rm P}_{\bf Y}$~probability content $\tau$ irrespective of~${\rm P}_{\bf Y}$ 
%
%The resulting quantile functions are supposed to  the images by  ${ Q}_{{\bf Y}}$ of 
%$${\mathbb C}_{\bf Y}(\tau ) =  \big\{
%{\bf y}=(y_1,\ldots,y_d) \vert \, {\rm P}_{{\bf Y}\text{\rm cop}}\big[F_1(y_1)\times\ldots\times F_d(y_d)
%\big]	\leq \tau \big\}\quad \tau\in (0,1),
%$$
%the images by 
%
and, for the same reason (lack of distribution-freeness) as above, the traditional definition of multivariate distribution functions does not allow for a satisfactory definition of quantiles and quantile regions, hence does not qualify as a valid multivariate extension of the univariate concepts. 
\subsection{The Rosenblatt transformation}\label{SecRosenblatt} 
Although the words {\it distribution} and {\it quantile functions} are not used in his 1952 note,\footnote{\cite{Rosenblatt52}, actually, is interested in multivariate versions of the Kolmogorov-Smirnov and von Mises Goodness-of-Fit tests.} the transformation proposed by \cite{Rosenblatt52}  and its inverse actually possess most of the attributes of  distribution and quantile functions. 

The so-called Rosenblatt transformation is a mapping 
$${\bf F}_{\text{$\bf Y$ {\rm Rosenblatt}}}: {\bf y}=(y_1,\ldots,y_d)\mapsto {\bf F}_{\text{ {\rm Rosenblatt}}}({\bf y})\coloneqq
\big(F_{\text{$\bf Y$ {\rm Rosenblatt}, 1}}({\bf y}),\ldots , F_{\text{$\bf Y${\rm Rosenblatt}, $d$}}({\bf y})
\big)$$
 from ${\mathbb R}^d$ to the unit cube $[0,1]^d$ pushing the distribution ${\rm P}_{\bf Y}\in{\mathcal P}$ of a $d$-dimensional random vector~${\bf Y}=(Y_1,\ldots,Y_d)$ forward to the uniform ${\rm U}_{[0,1]^d}$ over the unit cube. Denoting by $F_{Y_1}$ the traditional distribution function of $Y_1$  and by $F_{Y_j\vert Y_1,\ldots,Y_{j-1}}(\  \cdot \  \vert\, y_1,\ldots, y_{j-1})$ the distribution function of $Y_j$ conditional on $Y_1=y_1,\ldots,Y_{j-1}=y_{j-1}$, ${\bf F}_{\text{$\bf Y$ {\rm Rosenblatt}}}$ is defined recursively:%, with  
 \begin{align*}
{F}_{\text{{\rm Rosenblatt} }{\bf Y},\, 1 }({\bf y})\coloneqq& F_{Y_1}(y_1)\hspace{50mm} \\ 
{F}_{\text{{\rm Rosenblatt} }{\bf Y},\, 2 }({\bf y})\coloneqq& F_{Y_2\vert Y_1}(y_2\vert y_1) \\ 
\vdots\hspace{15mm}&\hspace{1mm}\vdots  \\
{F}_{\text{{\rm Rosenblatt} }{\bf Y},\, d }({\bf y})\coloneqq& F_{Y_d\vert Y_1,\ldots, Y_{d-1}}(y_d\vert\, y_1, \ldots, y_{d-1}).
\end{align*}

 If the density $f_{\bf Y}$ of ${\rm P}_{\bf Y}$ is continuous, ${\bf F}_{\text{$\bf Y$ {\rm Rosenblatt}}}$ is a homeomorphism from ${\mathbb R}^d$ to~$[0,1]^d$, and admits a continuous inverse~${\bf Q}_{\text{$\bf Y$ {\rm Rosenblatt}}}\coloneqq {\bf F}^{-1}_{\text{$\bf Y$ {\rm Rosenblatt}}}$ which pushes the uniform~${\rm U}_{[0,1]^d}$ forward to ${\rm P}_{\bf Y}$. Clearly,~${\bf F}_{\text{$\bf Y$ {\rm Rosenblatt}}}$ and ${\bf Q}_{\text{$\bf Y$ {\rm Rosenblatt}}}$ both characterize ${\rm P}_{\bf Y}$ and,  in dimension~$d=1$,  coincide with the traditional univariate concepts.  For~$d\geq~\!2$, they are not, in general, gradients of convex functions, hence violate (b) and (c);    ${\bf Q}_{\,\text{$\bf Y${\rm Rosenblatt}}}$, moreover, is highly non-equivariant under orthogonal transformations; ${\bf F}_{\text{$\bf Y$ {\rm Rosenblatt}}}$ and~${\bf Q}_{\,\text{ $\bf Y$ {\rm Rosenblatt}}}$, thus,    satisfy the desired properties (a) and (d), but not (b), (c), and (e).
 %satisfy the desired properties (a)--(d).   
 
 An additional weak point %, however,
  is the existence of $d!$ distinct versions (associated with the $d!$ possible orderings of $\bf Y$'s components), typically leading to markedly different solutions (no equivariance under orthogonal transformations). Moreover, when quantile regions of the form ${\bf Q}_{\text{$\bf Y$ {\rm Rosenblatt}}}({\scriptstyle {\mathbb C}}(\tau))$ are to be constructed, one runs into the problem of choo\-sing the adequate collection of nested regions~${\scriptstyle {\mathbb C}}(\tau)$,  $\tau\in (0,1)$  of $[0,1]^d$ independent of ${\rm P}_{\bf Y}$ and  such that~${\rm U}_{[0,1]^d}\big(
{{\ensuremath{\mathsmaller{\mathbb{C}}}}}(\tau) \big)=\tau$ for all~$\tau\in (0,1)$.   Such a choice is by no means obvious:  the ``south-east regions''  
\[ {\scriptstyle {\mathbb C}}(\tau)
%\text{$${{\ensuremath{\mathsmaller{\mathbb{C}}}}}$
%$} 
\coloneqq \{{\bf u}=(u_1,\ldots,u_d)\vert\, 0<u_j\leq \tau^{1/d},\ j=1,\ldots, d\},\quad \tau\in (0,1)
\]
(which, however, fail to satisfy the orthogonal equivariance property (e)), the  ``central cubes'' 
$$ {\scriptstyle {\mathbb C}}(\tau) \coloneqq \{{\bf u}=(u_1,\ldots,u_d)\vert\, (1-\tau^{1/d})/2<u_j\leq (1+\tau^{1/d})/2,\ j=1,\ldots, d\},\quad \tau\in (0,1)
$$ 
or the depth regions (see Section~\ref{Secdepth})
\begin{equation}\label{cubicdepth}
{\scriptstyle {\mathbb C}}(\tau) \coloneqq \{{\bf u}\vert\, D({\bf u}; {\rm U}_{[0,1]^d})\leq d(\tau)\},\quad \tau\in (0,1)
\end{equation}
where $D({\bf u}; {\rm U}_{[0,1]^d})$ denotes the Tukey depth\footnote{We concentrate here on halfspace or Tukey depth  as proposed in \citep{Tukey75}. However, in the wake of Tukey, several other depth functions have been proposed, such as  {\it simplicial depth} \citep{Liu90},   {\it zonoid depth} \citep{Mosler02}, or {\it projection depth} \citep{Zuo03}, to quote only a few; see also \cite{ZuoS00}.}  of $\bf u$ with respect to ${\rm U}_{[0,1]^d}$  and $d(\tau)$ is such that ${\rm U}_{[0,1]^d}\big(D({\bf u}; {\rm U}_{[0,1]^d})\leq d(\tau)
\big) = \tau$, are equally plausible  possibilities and, combined with the variety of possible Rosenblatt transformations, make practical interpretation challenging. This issue, which is common to all concepts of quantile functions with domain the unit cube,\footnote{This is not an  issue, however, when (as in spatial or center-outward quantiles) the unit ball, rather than the unit cube, is the quantile function domain: see Sections~\ref{Secspatial} and ~\ref{SecOT}.}  and the lack of orthogonal equivariance may explain why Rosenblatt's concept has only met with limited success in practice.  

\color{red}

\color{black}

\subsection{Elliptical or Mahalanobis quantiles}\label{EllSec} If $\bf Y\sim{\rm P}_{\bf Y}$ is elliptical, with center $\boldsymbol\mu$, positive definite scatter matrix $\boldsymbol\Sigma$, and radial density $f_{\text{rad}} $ nonvanishing on its support, the distribution~of 
$${\bf F}_{{\text{\bf Y}\, {\rm ell}}} ({\bf Y})\coloneqq F_{\text{rad}}\Big( (({\bf Y}- {\boldsymbol\mu})^\top{\boldsymbol\Sigma}^{-1}({\bf Y}- {\boldsymbol\mu}))^{1/2}\Big) \frac{{\boldsymbol\Sigma}^{-1/2}({\bf Y}- {\boldsymbol\mu})}{\big(({\bf Y}- {\boldsymbol\mu})^\top{\boldsymbol\Sigma}^{-1}({\bf Y}- {\boldsymbol\mu})\big)^{1/2}}$$
where $F_{\text{rad}}$ denotes the traditional cumulative distribution function associated with $f_{\text{rad}} $ and~${\boldsymbol\Sigma}^{1/2}$ is the symmetric root of $\boldsymbol\Sigma$,  is spherical uniform: ${\bf F}_{{\text{\bf Y}\, {\rm ell}}} ({\bf Y})\sim{\rm U}_d$. The mapping~${\bf y}\mapsto{\bf F}_{{\text{\bf Y}\, {\rm ell}}} ({\bf y})$, therefore, is a multivariate probability integral transformation pushing~${\rm P}_{\bf Y}$ forward to ${\rm U}_d$, hence an extended distribution function candidate: call it the {\it elliptical} or {\it Mahalanobis distribution function of} ${\rm P}_{\bf Y}$.  
Call {\it elliptical} or {\it Mahalanobis quantile function of}~${\rm P}_{\bf Y}$ its inverse ${\bf Q}_{{\text{\bf Y}\, {\rm ell}}}\coloneqq {\bf F}^{-1}_{{\text{\bf Y}\, {\rm ell}}}$ mapping ${\bf u}\in {\mathbb S}_d$ to ${\bf Q}_{{\text{\bf Y}\, {\rm ell}}}({\bf u})\coloneqq 
{\boldsymbol\Sigma}^{1/2}Q_{\text{rad}}(\Vert{\bf u}\Vert)\frac{\bf u}{\Vert{\bf u}\Vert}$ (where~$Q_{\text{rad}}$ stands for the traditional quantile function associated with  $f_{\text{rad}} $): ${\bf Q}_{{\text{\bf Y}\, {\rm ell}}}$ is pushing~${\rm U}_d$ forward to ${\rm P}_{\bf Y}$. The corresponding quantile  contours of order $\tau$
$${\mathcal C}_{{\text{\bf Y}\, {\rm ell}}}(\tau)\coloneqq {\bf Q}_{{\text{\bf Y}\, {\rm ell}}}(\tau{\mathcal S}_d)$$
are nested ellipsoids  centered at $\boldsymbol\mu$, with shape matrix ${\boldsymbol\Sigma}$, enclosing  quantile  regions ${\mathbb C}_{{\text{\bf Y}\, {\rm ell}}}(\tau)$ of order $\tau$. 

Under the assumption of an elliptical ${\rm P}_{\bf Y}$, ${\bf F}_{{\text{\bf Y}\, {\rm ell}}}$ and ${\bf Q}_{{\text{\bf Y}\, {\rm ell}}}$ clearly enjoy all the properties of distribution and quantile concepts, plus the affine-equivariance of ${\bf Q}_{{\text{\bf Y}\, {\rm ell}}}$.  These concepts underlie the rank-based techniques developed in \cite{HPd02, HPd02b, HPd04, HPd06}, \cite{HPd06b}, and \cite{HPd09}---the validity of which is limited to elliptical distributions.

 \cite{HLUSIM13},  \cite{HLUSIM15}, and \cite{HSiman16} propose to relax that limitation and extend the application of elliptical quantiles to general~${\rm P}_{\bf Y}\in{\mathcal P}$ with finite second-order moments. The question, however, arises of the choice of $\boldsymbol\Sigma$ and $Q_{\text{rad}}$, which, away from ellipticity, have little meaning. The idea, in   \cite{HLUSIM13}, is to replace them with the solution of an optimization problem. More precisely, denoting by $\rho_{\tau}$ the usual check function \eqref{rhoquant}, they define the elliptical quantile contour of order $\tau$ as the ellipsoid
\[
{\mathcal C}_{\text{$\bf Y$ \rm H\v{S}}} (\tau)\coloneqq 
\Big\{{\bf y}\in{\mathbb R}^d
\vert\,  {\bf y}^\top {\bf A}_{\tau} {\bf y} + {\bf y}^\top {\bf b}_{\tau} - c_{\tau} = 0
\Big\}
\]
where ${\bf A}_{\tau}\in {\mathbb R}^d\times {\mathbb R}^d$, ${\bf b}_{\tau}\in  {\mathbb R}^d\times {\mathbb R}$, and $c_{\tau}\in{\mathbb R}^+$ minimize, subject to ${\bf A}$ being symmetric and positive semidefinite with
determinant one (${\bf A}$ is thus a {\it shape matrix} in the sense of \cite{Paindaveine08}), the objective function
$$\Psi_{\tau}
 ({\bf A}, {\bf b}, c)\coloneqq {\rm E}\big[ \rho_{\tau} ({\bf Y}^\top{\bf AY} + {\bf Y}^\top{\bf b} - c)\big].
$$  
The  region ${\mathbb C}_{\text{$\bf Y$ \rm H\v{S}}} (\tau)$  enclosed by ${\mathcal C}_{\text{$\bf Y$ \rm H\v{S}}} (\tau)$ is the corresponding quantile region of order $\tau$.\smallskip 

%with the usual check function ?? (x) := x(? ? I(x < 0)) = max{(? ? 1)x, ? x}. 
The positive semidefiniteness of $\bf A$ 
and the condition on its determinant ensure that~${\mathcal C}_{\text{$\bf Y$ \rm H\v{S}}} (\tau)$ is %indeed 
an ellipsoid, centered at ${\boldsymbol \mu}_{\tau}\coloneqq - {\bf A}_{\tau}{\bf b}^{-1}_{\tau}/2$, with equation   
$({\bf y} - {\boldsymbol \mu}_{\tau} )^\top {\bf A}_{\tau} ({\bf y} -{\boldsymbol \mu}_{\tau}) = \kappa (\tau)$, where $\kappa (\tau) \coloneqq  c_{\tau} + {\bf b}_{\tau}^\top {\bf A}^{-1}_{\tau}
{\bf b}_{\tau} /4$. The condition det$({\bf A}_{\tau}) = 1$ can be viewed as an identification
constraint: for any $K > 0$, indeed, the triples $({\bf A}_{\tau}, {\bf b}_{\tau}, c_{\tau})$ and $(K{\bf A}_{\tau}, K{\bf b}_{\tau}, Kc_{\tau})$  define the same ellipsoid.

Under elliptical ${\rm P}_{\bf Y}$, the 
 contours  ${\mathcal C}_{\text{$\bf Y$ \rm H\v{S}}} (\tau)$ are nested ellipsoids, with     center~$\boldsymbol\mu$ and a matrix~${\bf A}$  that do not depend on $\tau$, coinciding with %${\mathbb C}_{\text{\rm ell} {\bf Y}} (\tau)$ and 
  ${\mathcal C}_{\text{$\bf Y$ \rm H\v{S}}} (\tau)$, and enjoying the same great properties. Under general ${\rm P}_{\bf Y}$,  the main properties of  ${\mathcal C}_{\text{$\bf Y$ \rm H\v{S}}} (\tau)$ and ${\mathbb C}_{\text{$\bf Y$ \rm H\v{S}}} (\tau)$ are summarized in \cite{HSiman16} (page 234): 
   in particular,  the quantile regions ${\mathbb C}_{\text{$\bf Y$ \rm H\v{S}}} (\tau)$ and contours ${\mathcal C}_{\text{$\bf Y$ \rm H\v{S}}} (\tau)$ are affine-equivariant and their ${\rm P}_{\bf Y}$ probability content is $\tau$ irrespective of~${\rm P}_{\bf Y}$. These are very strong properties for given $\tau$ which, however,   come at a cost: the ellipsoids~${\mathcal C}_{\text{$\bf Y$ \rm H\v{S}}} (\tau)$, typically, are no longer nested as $\tau$ increases---their centers and principal directions depend on $\tau$. Whether they characterize ${\rm P}_{\bf Y}$ is not known, but is quite unlikely. 
 
 \color{black}
 
\subsection{Depth and depth-based concepts}\label{Secdepth} The conctept of halfspace depth$^6$ (see, e.g., \cite{Nagy25}) was first introduced by \cite{Tukey75}, in a context of outlyingness detection. Depth and  quantile functions both pursue the same goal---namely, to order the points of a space in the absence of a canonical order---and the success of depth theories owes much to the lack of consensus regarding the   definitions of multivariate quantiles available   in the literature.   
%The interpretation of depth functions, their inverses,  and depth regions/contours in terms of distribution  and quantile functions, quantile contours, and quantile regions/ is more recent. 

The depth of a point ${\bf y}\in{\mathbb R}^d$ with respect to ${\rm P}_{\bf Y}$ is defined (see, e.g., \cite{Liu99})  as 
$$D({\bf y}; {\rm P}_{\bf Y})\coloneqq \inf_{{\bf u}\in{\mathcal S}_{d-1}} {\rm P}_{\bf Y}\big({\bf u}^\prime
{\bf Y}\leq {\bf u}^\prime
{\bf y} \big).
$$
In dimension $d=1$, it is easy to see that $F_{Y\pm}(y) = \big(2D(y;{\rm P}_Y )-1\big))\text{sign}(\frac{1}{2}-F_Y(y))$, establishing a direct relation between the depth function $y\mapsto D(y;{\rm P}_Y )$ and the center-outward distribution function $y\mapsto F_{Y\pm}(y)$. For general dimension $d$, the depth function mapping~${\bf y}\in{\mathbb R}^d$ to $D({\bf y}; {\rm P}_{\bf Y})$ ranges over $[1/2, 0)$, hence $1-2D({\bf y};{\rm P}_{\bf Y} )$ takes values in~$[0, 1)$.  Interpreting $1-2D({\bf y};{\rm P}_{\bf Y} )\eqqcolon F_{ {\bf Y} \text{\rm depth}}({\bf y})$ as a distribution function and its inverse~$Q_{ {\bf Y} \text{\rm depth}}\coloneqq F^{-1}_{ {\bf Y} \text{\rm depth}}$ as a quantile function, %yielding, with ${{\scriptstyle{\mathbb C}}}(\tau)=[0,\tau]$,
 it is tempting to adopt   the {\it depth regions}  
\[{\mathbb C}_{{\bf Y} \text{\rm depth}}(\tau)\coloneqq    \{{\bf y}\in{\mathbb R}^d \vert\,  D({\bf y};{\rm P}_{\bf Y} )\geq (1-\tau)/2 \} = Q_{ {\bf Y} \text{\rm depth}}([0, \tau]) \quad \tau\in [0, 1)
\]
and {\it depth contours}
\[{\mathcal C}_{{\bf Y}\text{\rm depth}}(\tau)\coloneqq  \{{\bf y}\in{\mathbb R}^d \vert\,  D({\bf y};{\rm P}_{\bf Y} ) = (1-\tau)/2 \} 
= Q_{ {\bf Y} \text{\rm depth}}(\tau)  \quad \tau\in [0, 1)
\]
as central quantile regions and contours of order $\tau$ (associated with ${{\scriptstyle{\mathbb C}}}(\tau)=[0,\tau]$).  However, in dimension two and higher, some of the essential properties of distribution and quantile functions are violated:   
\begin{enumerate}
\item[(a')]  $D(\,\cdot\,;{\rm P}_{\bf Y})$ does not characterize ${\rm P}_{\bf Y}\in{\mathcal P}$;\footnote{This important negative  result was obtained by \cite{Nagy21}; see also \cite{Nagy20}.} nor do $F_{ {\bf Y} \text{\rm depth}}$ and $Q_{ {\bf Y} \text{\rm depth}}$;
\item[(b')]   $D({\bf Y};{\rm P}_{\bf Y})$ is not uniform over $[0,1)$ and not even distribution-free; therefore, the proba\-bility content of ${\mathbb C}_{{\bf Y} \text{\rm depth}}(\tau)$, in general, is not  $\tau$, and strongly depends on ${\rm P}_{\bf Y}$.
% \item[--] 
% \item[--] 
% \item[--] 
 \end{enumerate}
 Some authors suggest to remedy (b') by relabeling the regions ${\mathbb C}_{{\bf Y} \text{\rm depth}}(\tau)$ and the corresponding contours ${\mathcal C}_{{\bf Y} \text{\rm depth}}(\tau)$ with their actual ${\rm P}_{\bf Y}$~probability contents, which is a bit ad hoc. But the violation  (a') of (a) is severely and  irrevocably redhibitory. While of independent interest, this concept of depth thus  fails to provide a satisfactory definition of  multivariate quantiles.

\subsection{Distance-based concepts: geometric or spatial,  $\rho$-, and Oja quantiles}\label{Secspatial}
This class of multivariate quantiles includes various very ingenious and  elegant   concepts based on extensions of the L$_1$ characterization (iii) of the univariate concept. We refer to  \cite{Oja99, Mard99, OjaRan04, Serfling08},  and the monograph by  \cite{Oja10} for  details, systematic  reviews and exhaustive lists of references, to \cite{Konen22} and \cite{DPdNagy}  for the more recent $\rho$- and Oja quantiles.  

The most popular representative  in this class  is  probably the concept of {\it spatial} or {\it geometric quantile} proposed by \cite{Chau96}, refined by \cite{Chakra01}  and extended to infinite-dimensional spaces by \cite{Chakra14}, to L$_p$ generalized check functions by \cite{Konen22}, and to  the hypersphere by \cite{Konen23}. 

Chaudhuri's proposal\footnote{An earlier version can be found in an unpublished  manuscript by \cite{DudleyKoltchinski92}. } consists of extending the univariate characterization (iii) to ${\mathbb R}^d$  by replacing, in  %the univariate version  
\eqref{Chaudhu1}, %of the  expected check function  characterization of center-outward contours, 
 $\tau v%\{\tau,\, -\tau\}=
 \in \tau{\mathcal S}_0$  with $\tau{\bf v}\in  \tau{\mathcal S}_{d-1}$, yielding the spatial quantile\footnote{Here  again, the assumption of  finite first-order moments for $\bf Y$ can be avoided by subtracting, in the expectation, $\Vert {\bf Y}\Vert + \tau{\bf v}^\prime{\bf Y} \tau $, which has no impact on the  $\arg\! \min_{{\bf y}\in{\mathbb R}^d}$.
}\begin{equation}\label{Chaudhu2}%
{\bf Q}_{{\bf Y}\text{\rm spatial}}(\tau{\bf v}) \coloneqq \arg\! \min_{{\bf y}\in{\mathbb R}^d}{\rm E}\big[
\Vert {\bf Y}-{\bf y}
\Vert +  \tau {\bf v}^\prime({\bf Y}- {\bf y})\big],\quad {\bf v}\in {\mathcal S}_{d-1}.
%\big]\label{Chaudhu1}
% {\rm E}\big[\rho_{(1\pm \tau)/2} (Y) = \arg\! \min_{y\in{\mathbb R}} {\rm E}\big[ \vert Y-y\vert  \pm\tau (Y-y)\big].
\end{equation}
Spatial quantiles, thus, are indexed by ${\bf u}=\tau{\bf v} \in{\mathbb S}_d$, which decomposes into\vspace{0mm}  ${\bf v} \!\coloneqq\!  \dfrac{\bf u}{\Vert\bf u\Vert} \in~\!{\mathcal S}_{d-1}$ (a direction, playing the role of a multivariate sign) and\vspace{0mm} $\tau\!\coloneqq~\!\!\Vert{\bf u}\Vert \in~\![0,1)$ (a probability). For~$d=1$, the spatial quantiles ${\bf Q}_{{\bf Y}\text{\rm spatial}}(\tau{\bf v})$, ${\bf v}\in{\mathcal S}_{d-1}$ reduce to  the center-outward quantiles $Q_{Y \pm}(\tau v)$,~$v\in~\{\pm 1\}$.

For ${\rm P}_{\bf Y}\in{\mathcal P}$, the spatial quantile functions ${\mathbb S}_d\ni{\bf u}\mapsto {\bf Q}_{{\bf Y}\text{\rm spatial}}({\bf u})\in{\mathbb R}^d$   are homeo\-morphisms, hence are continuously invertible,   characterizing  spatial distribution functions~${\mathbb R}^d\ni{\bf y}\mapsto  {\bf F}_{{\bf Y}\text{\rm spatial}}({\bf y})\coloneqq  {\bf Q}^{-1}_{{\bf Y}\text{\rm spatial}}({\bf y})$.  Actually, for non-atomic  ${\rm P}_{\bf Y}$, it holds  (see \cite{DudleyKoltchinski92, Kol93}, or Section~3 of \cite{Chau96}) that ${\bf Q}_{{\bf Y}\text{\rm spatial}}({\bf u})$ is the unique solution of ${\rm E}\big[( {\bf Y}-{\bf y})/\Vert  {\bf Y}-{\bf y}\Vert \big]= -{\bf u}
$ so that
$$
\big\Vert {\rm E}\big[
( {\bf Y}-
{\bf Q}_{{\bf Y}\text{\rm spatial}}({\bf u}))/
\Vert  {\bf Y}-{{\bf Q}_{{\bf Y}\text{\rm spatial}}({\bf u})}\Vert 
\big]\big\Vert 
=\Vert{\bf u}\big\Vert (= \tau \text{ for }{\bf u}=\tau{\bf v}).$$ 

As distribution and quantile functions, ${\bf F}_{{\bf Y}\text{\rm spatial}}$ and ${\bf Q}_{{\bf Y}\text{\rm spatial}}$ enjoy many  of the properties expected from such concepts. In particular, (a) they characterize ${\rm P}_{\bf Y}$ (Theorems 2.5 and 2.9 in  \cite{Kol97}),  a property extended to separable infinite-dimensional Hilbert spaces in Theorem~1.1 of \cite{GonzalezKonen2026}\footnote{\cite{GonzalezKonen2026} also establish the   much stronger (and quite surprising) result that even the scalar-valued mapping  ${\bf y}\mapsto \big\Vert {\bf F}_{{\bf Y}\text{\rm spatial}}(\bf y)
\big\Vert$ characterizes ${\rm P}_{\bf Y}$.}
 and (e) yield equivariant (under shift, global rescaling, and orthogonal transformations) spatial quantile regions and contours 
$${\mathbb C}_{{\bf Y}\text{\rm spatial}}(\tau)\coloneqq {\bf Q}_{{\bf Y}\text{\rm spatial}}(\tau{\mathbb S}_d)\quad\text{and}\quad {\mathcal C}_{{\bf Y}\text{\rm spatial}}(\tau)\coloneqq {\bf Q}_{{\bf Y}\text{\rm spatial}}(\tau{\mathcal S}_{d-1}), \quad \tau\in [0,1);$$ 
see \cite{Konen25}.

No matter how appealing these spatial concepts are, however,   ${\bf F}_{{\bf Y}\text{\rm spatial}}({\bf Y})$  fails to be distribution-free. Property (b2), thus, is violated, hence also (b) and the crucial pro\-perty~(d): the probability content of  ${\mathbb C}_{{\bf Y}\text{\rm spatial}}(\tau)$ (the spatial quantile region of order $\tau$) is not~$\tau$, and still generally depends on ${\rm P}_{\bf Y}$.  Reindexing the quantile regions ${\mathbb C}_{{\bf Y}\text{\rm spatial}}(\tau)$ with their probability contents ${\rm P}_{\bf Y}\big({\mathbb C}_{{\bf Y}\text{\rm spatial}}(\tau)
\big)$ is sometimes proposed as a remedy to this problem. Such relabeling, however, is quite ad hoc. Actually, it turns {\it any} collection of nested regions  partitioning~${\mathbb R}^d$,  no matter how arbitrary it may be,  into an alleged collection of quantile regions, thereby stripping the concept of all meaningful interpretation. 

Moreover, it has been shown \citep{GirardStupfler17, HallinKonen} that spatial quantiles may exhibit quite strange behaviors as the value of $\tau$ increases: the norm of extreme geometric quantiles is  largest in the
direction where the variability of $\bf Y$ is  smallest (see Figure~1 and~2 in \cite{HallinKonen}),  and upper quantiles   exit the support of ${\rm P}_{\bf Y}$ when the latter is bounded (see Theorem 2.1  in \cite{GirardStupfler17} and 
Theorems~2--3 in \cite{DPdVirta}). 

\subsection{Some other multivariate quantile concepts} A variety of other concepts have been proposed in the literature, which are more or less related to the previous ones. 

Combining directional projections and depth features, \cite{KongMiz08} and \cite{KongMiz12} construct quantile-like contours which they use for inference (insisting, however, that they should not be interpreted as quantile contours in the strict sense). The corresponding population contours are similar to the ones considered in \cite{HPS10}. \cite{Fraiman12} also consider quantiles of directional projections that cover the infinite-dimensional case. \cite{Kol97} proposes quantiles based on the inversion of ``surrogate distribution functions'' of the form $G_{F_{\bf Y}},\ {\mathbb R}^d\to {\mathbb R}^d$ which characterizes $F_{\bf Y}$ and such that  the mapping ${\bf y}\mapsto G_{F_{\bf Y}}({\bf y})$ has an inverse, whose va\-lues he interprets as multivariate quantiles for ${\rm P}_{\bf Y}$. This approach, in a sense, is a precursor of the measure-transportation approach developed in Section~\ref{SecOT} but does not impose on $G_{F_{\bf Y}}$ any of the conditions that make $G^{-1}_{F_{\bf Y}}$ a convincing concept of quantile function. 
\cite{AbdTh92} can be considered as a forerunning variant of \cite{Chau96}.  Finally, \cite{Hett92} propose a definition based on the minimization of various measures of deviation of $\bf Y$ from a centerpoint, and, in a somewhat different perspective, \cite{EinmahlMason92, Polonik97} construct very general ``generalized quantile processes.''
 We refer to \cite{Serfling08} for details.

 None of these concepts, however, satisfies properties (b), and (d), and most of them also fail to meet property (a). 

 More recently, \cite{DPdNagy} introduced, under the name {\it Oja quantiles}, an ingenious concept which is to the Oja median what spatial quantiles are to the spatial median. Oja quantiles in dimension $d=1$ reduce to the univariate center-outward quantiles and constitute an extension of the spatial concept.  The Oja quantile function ${\bf Q}_{\text{\rm NPd}\,{\bf Y}}$ is, under mild regularity assumptions,  a homeomorphism between the open unit ball and ${\mathbb R}^d$, hence admits an inverse~${\bf F}_{\text{\rm NPd}\,{\bf Y}}$---a distribution function (which Nagy and Paindveine call a {\it rank function}). The Oja quantile and distribution functions characterize ${\rm P}_{\bf Y}$ (property (a)). They  enjoy affine-invariance/equivariance and, under ellipticity, the Oja quantile contours coincide with the natural ellipsoidal density contours (hence the elliptical or Mahalanobis contours discussed in Section~\ref{EllSec}).  The authors also demonstrate how their concept surpasses the more traditional spatial concept in many respects. On the other hand, the definitions of~${\bf Q}_{\text{\rm NPd}\,{\bf Y}}$ and~${\bf F}_{\text{\rm NPd}\,{\bf Y}}$ require finite moments of order one and, for large $\tau$ values,   exhibit the same weird behavior as spatial quantiles; whether  ${\bf F}_{\text{\rm NPd}\,{\bf Y}}({\bf Y})$ is distribution-free (property (b2), guaranteeing the~${\rm P}_{\bf Y}$ probability content of the resulting quantile region of order $\tau$ to be $\tau$)  is unknown, but seems unlikely.

\section{Measure-transportation-based distribution and quantile functions}\label{SecOT}\!\!\!\!\!\!\footnote{We refer to  \cite{Vil03, Vil09} 
for   background reading and an authoritative account of measure transportation, along with the two volumes by \cite{RR98} and the monograph by \cite{PanaZ},  the scope of which is  closer to  statistical concerns.  
}
\subsection{Center-outward distribution and quantile functions}\label{Sec41}
The  multiple attempts described in Section~\ref{Sec3} to extend the univariate definition  of  quantile functions,  thus, all fail to satisfy the fundamental properties expected from the concept of quantile. %have met with limited success. 
Except for the Rosenblatt transformation (Section~\ref{SecRosenblatt}), which degenerates in dimension one, 
% and hence has no clear univariate counterpart,
   these attempts are based on multivariate versions of the univariate characterizations (iii), (iv), and~(v). Only characterizations (i) and (ii),  which define~$F_{Y\pm}$ and $Q_{Y\pm}$ as monotone transports to and from the uniform distribution ${\rm U}_1$ over the unit ball, remained unexplored until very recently.

%$$F_{Y\pm}\# {\rm P}_Y = {\rm U}_1\quad \text{ and } \quad Q_{Y\pm}\#  {\rm U}_1 ={\rm P}_Y ,
%$$ 
%thus, remained unexplored\footnote{Koltchinski}  \citep{chernozhukov2017monge} and \citep{hallin2021distribution}. 

This is all the more surprising since  characterizations (i)--(ii) extend almost {\it verbatim} from dimension one to dimension $d$. Indeed, a  famous theorem by \cite{McCann95} implies that, denoting by ${\rm U}_d$ the spherical\footnote{The spherical uniform distribution ${\rm U}_d$ over ${\mathbb S}_d$ is the distribution of ${\bf U}\coloneqq t{\bf V}$ where $\bf V$ and $t\coloneqq \Vert\bf U\Vert$ are independent, $\bf V$  is uniform over ${\mathcal S}_{d-1}$, and $t$ is  uniform over $[0,1]$. For $d=1$, it coincides with the Lebesgue uniform over ${\mathbb S}_1=(-1,1)$.} uniform over ${\mathbb S}_d$, there exists a unique\footnote{Actually, uniqueness only holds almost surely, that is, up to a set of $\bf Y$ values contained in a Lebesgue null set. For the sake of simplicity, we henceforth omit mentioning this.} gradient of convex function pushing ${\rm P}_{\bf Y}\in{\mathcal P}$ forward to  ${\rm U}_d$: denote it by ${\bf F}_{{\bf Y}\pm}$ and call it the {\it center-outward distribution function of}  ${\rm P}_{\bf Y}$. \cite{Fig18} shows that~${\bf F}_{{\bf Y}\pm}$ is a homeomorphism between the pointed open unit ball ${\mathbb S}_d\!\setminus\!\{{\boldsymbol 0}\}$ and its image by~${\bf F}_{{\bf Y}\pm}$, hence admits, on ${\mathbb S}_d\!\setminus\!\{{\boldsymbol 0}\}$, an inverse~${\bf F}^{-1}_{{\bf Y}\pm}\eqqcolon {\bf Q}_{{\bf Y}\pm}$: call  it %${\bf Q}_{{\bf Y}\pm}$
 the {\it center-outward quantile function of}  ${\rm P}_{\bf Y}$. 

Due to the invariance under orthogonal transformation of ${\rm U}_d$, the choice of the orthogonal-invariant collection ${\scriptstyle{\mathbb C}}(\tau)\coloneqq 
\tau \overline{\mathbb S}_d$ (the closed ball with radius $\tau\in (0,1)$, hence ${\rm U}_d$~probability content $\tau$) as quantile regions for~${\rm U}_d$ is a  natural one. This choice  provides {\it center-outward quantile regions and contours of order $\tau$} of the form 
$${\mathbb C}_{{\bf Y} \pm}(\tau)\coloneqq {\bf Q}_{{\bf Y}\pm}({\scriptstyle{\mathbb C}}(\tau)) \quad\text{ 
%$$
and, with ${\scriptstyle{\mathcal C}}(\tau)\coloneqq \tau{\mathcal S}_{d-1}$,}\quad  
%$$
{\mathcal C}_{{\bf Y} \pm}(\tau)\coloneqq {\bf Q}_{{\bf Y}\pm}({\scriptstyle{\mathcal C}}(\tau)) %= 
,$$
respectively, for $\tau\in (0,1)$.   If we include the origin in its range, ${\bf F}_{{\bf Y}\pm}$ is no longer a  homeomorphism; for $\tau = 0$, we therefore define ${\mathbb C}_{{\bf Y} \pm}(0)$ as  $\bigcap_{\tau>0}{\mathbb C}_{{\bf Y} \pm}(\tau )$, which is not necessarily a singleton: call it the {\it center-outward median set of} ${\rm P}_{\bf Y}$.   

This is the approach taken in\footnote{As mentioned before, \cite{Kol97} can be considered as an early precursor of measure-transportation-based methods in the area.}  \cite{chernozhukov2017monge}\footnote{In \cite{chernozhukov2017monge}, however, the uniform over the unit cube $[0,1]^d$ is privileged, raising the same issue in the choice of the  ${\scriptstyle{\mathbb C}}(\tau)$ regions as in \eqref{cubicdepth}. The population concepts adopted in  \cite{ghosal22} are the same as here and in \cite{hallin2021distribution}, under a slightly different terminology and more flexibility in the choice of the reference distribution. Their empirical treatment, however, differs from that of \cite{hallin2021distribution}.} and \cite{hallin2021distribution}, followed by \cite{ghosal22}$^{16}$ and \cite{Bercu24}.  \cite{FaugRusch17} also propose a similar concept, combined with a preliminary copula transform.  See \cite{HAnnRev2021} for a review.

Contrary to all previous concepts (several of which, such as  spatial  quantiles, also involve transports to and from the unit ball ${\mathbb S}_d$, albeit not to and from the spherical uniform ${\rm U}_d$), the measure-transportation-based center-outward distribution and quantile functions ${\bf F}_{{\bf Y}\pm}$ and~${\bf Q}_{{\bf Y}\pm}$ enjoy all the properties expected from these concepts: 
\begin{enumerate}
\item[(a)] by construction, ${\bf Q}_{{\bf Y}\pm}$ characterizes ${\rm P}_{\bf Y}$ and,  being ${\rm P}_{\bf Y}$'s inverse, so does ${\bf F}_{{\bf Y}\pm}$;
\item[(b)--(c)] still by construction, ${\bf F}_{{\bf Y}\pm}$ and ${\bf Q}_{{\bf Y}\pm}$, as gradients of convex functions, are (cyclically) monotone, and  ${\bf F}_{{\bf Y}\pm}({\bf Y})\sim {\rm U}_d$, hence is distribution-free, which entails that 
\begin{enumerate}
\item[(d)] the ${\rm P}_{\bf Y}$~probability content of ${\mathbb C}_{{\bf Y} \pm}(\tau)$ is $\tau$ for all $\tau\!\in\![0,1)$, irrespective of~${\rm P}_{\bf Y}\!\in~\!\!{\mathcal P}\!$, justifying the terminology quantile region ``of order $\tau$,''~and 
%\item[(f)] $\{({\bf Q}_{{\bf Z}\pm}\circ {\bf F}_{{\bf Y}\pm})\#{\rm P}_{\bf Y} \vert\, {\bf Z}\sim{\rm P}\in {\mathcal P}\} ={\mathcal P}$;
\end{enumerate}
\item[(e)] equivariance of quantile regions and quantile contours under shift, rescaling, and ortho\-gonal transformations: for any  $\tau\in (0,1)$, any ${\bf a}\in{\mathbb R}^d$, any $b\in {\mathbb R}\setminus\{0\}$, and any $d\times d$ orthogonal matrix ${\bf O}$, ${\mathbb C}_{{\bf a}+b{\bf OY} \pm}(\tau) = {\bf a}+b{\bf O}{\mathbb C}_{{\bf Y} \pm}(\tau)$ and ${\mathcal C}_{{\bf a}+b{\bf OY} \pm}(\tau) = {\bf a}+b{\bf O}{\mathcal C}_{{\bf Y} \pm}(\tau)$. 
%C()(\tau)={\bf a+b{\bf OY}(\tau) and C(a+bY)±(?)=a+bCY±(?);
%\item[(f)] 
\end{enumerate}

Measure-transportation-based quantiles, thus, may not be the only solution extending to~${\mathbb R}^d\ (d>1)$ the essential properties 
%(a)--(c)
 of the traditional univariate  concept, but they are the only ones so far. 

%We refer to  \cite{Vil03, Vil09} 
%for   background reading and an authoritative account of measure transportation, along with the two volumes by \cite{RR98} and the monograph by \cite{PanaZ},  the scope of which is  closer to  statistical concerns.  
%
%

\subsection{Quantiles on manifolds} Multidimensional real spaces are not the only ones where, due to the absence of a canonical ordering, no obvious extension of the concepts of distribution and quantile functions are available. Here again, measure transportation ideas yield  concepts that enjoy the essential properties (a)--(c) of their univariate counterparts.  The very important special case of hyperspheres (directional data)   has been considered in \cite{HLV24}, more general Riemannian manifolds in \cite{HallinLiuRiemann24}. 

\subsubsection{Quantiles on the hypersphere ${\mathcal S}_{d-1}$}\label{Sec421}

 Inference for directional data has a long tradition. While the circle (the 1-sphere) admits a canonical {\it circular} ordering, nothing comparable is available on the hypersphere ${\mathcal S}_{d-1}\coloneqq\{{\bf v}\vert\, \Vert{\bf v}\Vert\ = 1\}$ (the $(d-1)$-sphere   in ${\mathbb R}^d$), $d\geq 3$). As a result, inference about directional data is often limited to rotationally invariant distributions. Although extensions of rotational symmetry yielding, after projection onto  tangent spaces, elliptically symmetric rather than spherical distributions have been proposed in \cite{Kent82}, \cite{Scealy20}, or \cite{Garcia20}, 
%---the directional counterpart of spherical symmetry---despite empirical evidence showing that many datasets are not rotationally symmetric.
this places substantial limits on the validity of the quantile concepts proposed, e.g., by \cite{Ley14}, so we will not discuss them any further. The spatial quantiles on ${\mathcal S}_{d-1}$ developed by  \cite{Konen23} have been discussed in Section~\ref{Secspatial}.

Let ${\mathcal P}_{{\mathcal S}_{d-1}}$ denote the class of distributions on ${{\mathcal S}_{d-1}}$    
%equipped with the $\sigma$-field associated with its Riemannian metric
 which are absolutely continuous with respect to the surface measure $\lambda_{{\mathcal S}_{d-1}}$ on ${\mathcal S}_{d-1}$ equipped with the geodesic distance, with  density (with respect to $\lambda_{{\mathcal S}_{d-1}}$) bounded from above and away from zero. Based on measure transportation results for hyperspheres, 
% \footnote{These results are summarized in Propositions~1 and~2 of \cite{HLV24}, where we refer for details. } 
 \cite{HLV24} propose to call distri\-bution and quantile functions of ${\rm P}_{\bf Y}\in{\mathcal P}_{{\mathcal S}_{d-1}}$ the unique$^4$ mapping ${\bf F}_{{\mathcal S}_{d-1}{\bf Y}}$ from~${\mathcal S}_{d-1}$~to~${\mathcal S}_{d-1}$ pushing ${\rm P}_{\bf Y}$ forward to the uniform ${\rm U}_{{\mathcal S}_{d-1}}$ over ${\mathcal S}_{d-1}$ 
 and its inverse\footnote{It follows from Propositions~1 and~2 in \cite{HLV24}  that ${\bf F}_{{\mathcal S}_{d-1}{\bf Y}}$ exists, is unique$,^4$ and a homeomorphism.}~${\bf Q}_{{\mathcal S}_{d-1}{\bf Y}}\coloneqq {\bf F}^{-1}_{{\mathcal S}_{d-1}{\bf Y}}$ pushing ${\rm U}_{{\mathcal S}_{d-1}}$ forward to ${\rm P}_{\bf Y}$. 
 
 The construction of quantile regions and contours, as in ${\mathbb R}^d$, further requires the choice of a collection of nested regions $\scriptstyle{\mathbb C}(\tau)$ with ${\rm U}_{{\mathcal S}_{d-1}}$~probability content $\tau$, $\tau\in [0,1]$.\footnote{Since ${\mathcal S}_{d-1}$ is compact, the range for $\tau$ can include $\tau = 1$.} Sensible choices of these regions (which are quantile regions for  ${\rm U}_{{\mathcal S}_{d-1}}$) should be based on the symmetry properties of ${\rm U}_{{\mathcal S}_{d-1}}$, which is invariant under the group ${\mathcal G}_{\text{orth}}$ of orthogonal transformations of ${\mathbb R}^d$. Invariance with respect to the full group ${\mathcal G}_{\text{orth}}$, however, is impossible 
% : $\scriptstyle{\mathbb C}(\tau)$, thus, should be invariant under the same group, which is impossible
  unless $\tau =0$ or $1$.   \cite{HLV24} therefore propose to concentrate on subgroups and, fixing a  pole  $\thetab$ such as an element of the Fr\' echet mean set of ${\rm P}_{\bf Y}$,  consider the  subgroup~${\cal G}_{\thetab}\subset {\mathcal G}_{\text{orth}}$ of rotations and symmetries with axis~${\bf F}_{{\mathcal S}_{d-1}{\bf Y}}(\thetab)$. This pole, in turn, characterizes an {\it equator}~$\text{Eq}_{\thetab}\coloneqq\{{\bf v}\in{\mathcal S}_{d-1}\vert\, {\bf v^T}{\bf F}_{{\mathcal S}_{d-1}{\bf Y}}(\thetab)  = 0\}$ {(a $(d-2)$-sphere)},   
{\it spherical caps} centered at~$ {\bf F}_{{\mathcal S}_{d-1}{\bf Y}}(\thetab) $,  
%with~$\lambda_{{\mathcal S}_{d-1}}(\scriptstyle{\mathbb C}(\tau)) = 2\pi^{d/2}\tau/\Gamma(d/2)$
 and {\it equatorial strips} symmetric with respect to $\text{Eq}_{\thetab}$, all of  which are invariant under ${\cal G}_{\thetab}$. 
 
 The spherical cap with ${\rm U}_{{\mathcal S}_{d-1}}$ probability content $\tau$ is 
  \[
{\scriptstyle{\mathbb C}}_{\text{cap}}(\tau) \coloneqq 
\left\{{\bf v}\in {\cal S}^{d-1} \vert\, F_{*}\left({\bf v}^\top {\bf F}_{{\mathcal S}_{d-1}{\bf Y}}(\thetab) \right)\geq 1-\tau
\right\} \quad 0\leq\tau \leq 1,\ 
\]
with boundary\vspace{-2mm}
\begin{equation*}\label{paraldef}
{\scriptstyle{\mathcal C}}_{\text{cap}}(\tau) \coloneqq 
\left\{{\bf v}\in {\cal S}^{d-1} \vert\, F_{*}\left({\bf v}^\top {\bf F}_{{\mathcal S}_{d-1}{\bf Y}}(\thetab) \right)= 1-\tau
\right\} \quad 0<\tau< 1, 
\end{equation*}
{(a $(d-2)$-sphere)} where 
\begin{equation*}\label{eq:unif}
u\mapsto F_{*}(u)\coloneqq  {\int_{-1}^u (1-s^2)^{(d-3)/2} \; \dd s}{\Big/}{\int_{-1}^1 (1-s^2)^{(d-3)/2} \;  \dd s}, \quad -1 \leq u \leq 1%\vspace{-1mm}
\end{equation*}
 is the distribution function of  
%the random variable 
${\bf V}^\top {\bf F}_{{\mathcal S}_{d-1}{\bf Y}}(\thetab) $ where~${\bf V}\sim {{\rm U}_{{\mathcal S}_{d-1}}}$. The resulting   quantile contour and quantile region of order $\tau$ of ${\bf Y} \sim {\rm P}_{\bf Y}$  then are 
$%\begin{equation}\label{qcontdef}
{\mathcal C}_{\text{cap}{\bf Y}}(\tau) 
\coloneqq {\bf Q}_{{\mathcal S}_{d-1}{\bf Y}}
\left({\scriptstyle{\mathcal C}}_{\text{cap}}(\tau)\right)
%=
% \left\{{\bf z}\in \mathcal{S}^{d-1} : F_{*} \left(({\bf F} ({\bf z})^\top {\bf F}(\thetab_{\rm M}))\right) =1-\tau
%\right\} )
$ %\vspace{-2mm}\end{equation}
and~${\mathbb C}_{\text{cap}{\bf Y}}(\tau) 
\coloneqq {\bf Q}_{{\mathcal S}_{d-1}{\bf Y}}
\left({\scriptstyle{\mathbb C}}_{\text{cap}}(\tau)\right)$, 
%=
% \left\{{\bf z}\in \mathcal{S}^{d-1} : F_{*} \left(({\bf F}  ({\bf z}))^\top {\bf F}(\thetab_{\rm M})\right) \geq 1- \tau
%\right\},\vspace{-2mm}\end{equation}
% by $\bf Q$ of ${\mathcal C}_\tau^{{\bf U}}$ and ${\mathbb C}_\tau^{{\bf U}}$,
  respectively. Since ${\bf Q}_{{\mathcal S}_{d-1}}$ is a measure-preserving transformation pushing ${\rm U}_{{\mathcal S}_{d-1}}$ forward to ${\rm P}_{\bf Y}$, the~${\rm P}_{\bf Y}$ probability content of ${\mathbb C}_{\text{cap}{\bf Y}}(\tau)$ is $\tau$, as expected. 
  
Note that  $F_*({\bf v}^\top {\bf F}_{{\mathcal S}_{d-1}{\bf Y}}(\thetab))$, ${\bf v}\in{\mathcal S}_{d-1}$ is  a measure of ${\bf v}$'s latitude, ranging from~0 (for~${\bf v} =-  {\bf F}_{{\mathcal S}_{d-1}{\bf Y}}(\thetab)$) to 1 (for ${\bf v} =  {\bf F}_{{\mathcal S}_{d-1}{\bf Y}}(\thetab)$) and scaled in such a way that the latitude of~${\bf V}\sim {{\rm U}_{{\mathcal S}_{d-1}}}$ is uniform over $[0,1]$. Therefore, for the hypersphere $ {\cal S}^{d-1}$ {equipped with the uniform distribution~$ {{\rm U}_{{\mathcal S}_{d-1}}}$} and the pole~$ {\bf F}_{{\mathcal S}_{d-1}{\bf Y}}(\thetab) $, the contour ${\scriptstyle{{\mathcal C}}}_{\text{cap}}(\tau)$ plays the role of a {\it parallel} of order $\tau$. Similarly, ${\mathcal C}_{\text{cap}{\bf Y}}(\tau)$ plays, for ${\cal S}^{d-1}$ now equipped with the  (non-uniform)  distribution~${\rm P}_{\bf Y}$ and the pole~$\thetab$, the role of a ``curvilinear'' parallel   of order $\tau$ adapted to~${\rm P}_{\bf Y}$  yielding, for $\tau=1/2$, the ``curvilinear equator''  ${\mathbb C}_{\text{cap}{\bf Y}}(1/2)= {\bf Q}_{{\mathcal S}_{d-1}{\bf Y}}(\text{Eq}_{\thetab})$. 

Rather than the spherical caps ${\scriptstyle{\mathbb C}}_{\text{cap}}(\tau)$ centered at  ${\bf F}_{{\mathcal S}_{d-1}{\bf Y}}(\thetab)$, one may prefer the equa\-torial strips 
\[{\scriptstyle{\mathbb C}}_{\text{Eq}}(\tau)\coloneqq 
\Big\{{\bf v}\in {\cal S}^{d-1} \Big\vert\,  \frac{1-\tau}{2}\leq F_{*}\left({\bf v}^\top {\bf F}_{{\mathcal S}_{d-1}{\bf Y}}(\thetab) \right)\leq \frac{1+\tau}{2}
\Big\} \quad 0\leq\tau \leq 1
\]
 comprised between ${\scriptstyle{\mathcal C}}_{\text{cap}}(
(1+\tau)/2)$ and~${\scriptstyle{\mathcal C}}_{\text{cap}}(
(1-\tau)/2)$ and  centered at $\text{Eq}_{\thetab}$,  with ${\rm U}_{{\mathcal S}_{d-1}}$ proba\-bility content $\tau$. 
The resulting quantile regions of order $\tau$ are ${\mathbb C}_{\text{Eq}{\bf Y}}(\tau)\coloneqq
 {\bf Q}_{{\mathcal S}_{d-1}{\bf Y}}({\scriptstyle{\mathbb C}}_{\text{Eq}}(\tau))$, with quantile contours 
\[{{\mathcal C}}_{\text{Eq}{\bf Y}}(\tau)\coloneqq {{\mathcal C}}_{\text{cap}{\bf Y}}((1+\tau)/2)\bigcup  {{\mathcal C}}_{\text{cap}{\bf Y}}((1-\tau)/2). 
\]
  
 The  essential properties (a) and (b), hence also (c) and (d)  of quantiles readily follow from the definitions of ${\bf F}_{{\mathcal S}_{d-1} {\bf Y}}$, ${\bf Q}_{{\mathcal S}_{d-1} {\bf Y}}$,   the quantile contours~${\mathcal C}_{\text{cap}{\bf Y}}(\tau)$, ${\mathcal C}_{\text{Eq}{\bf Y}}(\tau)$, and the quantile  regions~${\mathbb C}_{\text{cap}{\bf Y}}(\tau)$, ${\mathbb C}_{\text{Eq}{\bf Y}}(\tau)$   and the continuity of ${\bf F}_{{\mathcal S}_{d-1}}$ and~${\bf Q}_{{\mathcal S}_{d-1}}$; as for~(e), it takes the form of equivariance under orthogonal transformations, shifts and rescaling of the plunging space ${\mathbb R}^d$. 
 % details are left to the reader.
%\begin{proposition} \label{characterize}
%Let ${\bf Z}\sim{\rm P}^{\bf Z}  \in \mathfrak{P}_d$ have distribution and quantile functions~${\bf F}$ and~${\bf Q}$, respectively. Then, 
%\begin{compactenum}
%\item[(i)] ${\bf F}$ entirely characterizes ${\rm P}^{\bf Z}$,  ${\bf F}({\bf Z})\sim{\bf F}\#{\rm P}^{\bf Z}=~\!{\rm U}_{{\mathcal S}_{d-1}}$, and  $F_{*} \left(({\bf F} ({\bf Z}))^\top {\bf F}(\thetab_{\rm M})\right) \sim~\!{\rm U}_{[0,1]}$;
%\item[(ii)] ${\bf Q}$ entirely characterizes ${\rm P}^{\bf Z}$ and ${\bf Q}({\bf U})\sim{\bf Q}\#{\rm U}_{{\mathcal S}_{d-1}}=~\!{\rm P}^{\bf Z}$
%;% and ${\bf Q}({\bf Z})$ is independent of ${\rm P}^{\bf Z}$;
%\item[(iii)] the quantile contours $\mathcal C_\tau$, $\tau\in[0,1]$ are continuous;  the quantile regions~${\mathbb C}_\tau$ are   closed, connected, and nested; their intersec\-tion~$\bigcap_{\tau\in[0,1]} \mathbb{C}_\tau$ is the {\rm directional median}~$\thetab_{\rm M}$;
%\item[(iv)] the probability content ${\rm P}^{\bf Z}\left(\mathbb C_\tau
%\right)$ of ${\mathbb C}_\tau$, $\tau\in[0,1]$,  is $\tau$, irrespective of $\rm P ^{\bf Z}\in  \mathfrak{P}_d$. \vspace{-2mm}
%\end{compactenum}
%\end{proposition}

 \subsubsection{Quantiles on compact Riemannian manifolds}
 
 Increased attention has been given recently to the statistical analysis of variables with values on more complex Riemannian manifolds, such as the $p$-torus (a product of $p$ circles) or the polyspheres (the product of $k$ hyperspheres). The ideas developed for hyperspheres in the previous section apply to  such manifolds, with some important differences, though.

The hypersphere has the property that all of its closed connected contours\footnote{The word ``contour'' here is used in the loose sense of a  continuous curve that has no end point and encloses a connected region.} and, more particularly, the quantile contours ${\mathcal C}_{\text{cap}{\bf Y}}(\tau)$, are homotopic to a single point (such as $\thetab$). This is no longer the case  in general, and   closed contours of various homotopy types typically coexist. Required measure transportation results (e.g.,  the continuity of optimal transports), moreover, are not available for all types of manifolds---the most recent ones \citep{Rusch96, SCHACH08, McCann01} 
 are limited to connected $C^4$-smooth compact Riemannian manifolds without boundaries and with  everywhere nonnegative sectional
curvature. This includes the hypersphere, the $p$-torus, and the polyspheres but not the (poly)cylinders, Stiefel, or Grassmann manifolds.

As a simple example,  consider  the 2-torus ${\mathcal T}_2$, that is, the product ${\mathcal S}_1^1\times {\mathcal S}_1^2$ of two circles,  interchangeably called  {\it toroidal} and {\it poloidal}, equipped with its {\it Riemannian metric}~$g$. Continuous closed contours on ${\mathcal T}_2$ can be homotopic to one single point (thus defining {\it toroidal caps}), to the toroidal circle ${\mathcal S}_1^1$, or to the poloidal one ${\mathcal S}_1^2$ (characterizing {\it toroidal} and {\it poloidal equatorial strips} respectively). Each type yields distinct collections of quantile contours and regions, with complementary interpretations. 

Let ${\mathcal P}_{{\mathcal T}_2}$ denote the family of distributions on ${\mathcal T}_2$ that are absolutely continuous and admit a bounded (from above and from below) density with respect to the surface measure~$\lambda_{{\mathcal T}_2}$. The existence, uniqueness, and homeomorphic nature  of  the optimal\footnote{On compact manifolds, (cyclically) monotone transports and optimal (with respect to the squared geodesic distance) transports  coincide.} transport pushing~${\rm P}_{\bf Y}\in~\!{\mathcal P}_{{\mathcal T}_2}$ forward to the uniform ${\rm U}_{_{{\mathcal T}_2}}$ (a rescaled version of $\lambda_{{\mathcal T}_2}$ integrating to one) follow from the aforementioned results on measure transportation on manifolds. Call distribution function of ${\rm P}_{\bf Y}$ this optimal transport,  quantile function its inverse, and denote them as ${\bf F}_{{\mathcal T}_2 {\bf Y}}$ and ${\bf Q}_{{\mathcal T}_2 {\bf Y}}$, respectively: both are characterizing ${\rm P}_{\bf Y}$ (property (a)).

If quantile regions are to be constructed, however, this is not enough, and families of nested sets ${\scriptstyle{\mathbb C}}(\tau)$ with ${\rm U}_{_{{\mathcal T}_2}}$ probability content $\tau\in[0,1]$ have to be chosen.  Based on the invariance properties (with respect to toroidal symmetries and rotations, poloidal symmetries and rotations, or both), \cite{HallinLiuRiemann24} propose three choices. These three choices involve a centering---for instance,  an element $\thetab = (\theta_1,\theta_2)$, with $\theta_1\in {\mathcal S}_1^1$  and $\theta_2\in {\mathcal S}_1^2$, of the Fr\' echet mean set of~${\rm P}_{\bf Y}$. As a Fr\' echet mean, the point $\thetab$ can be considered as ``central'' in~${\rm P}_{\bf Y}$: call it a {\it pole}. Similarly, the curve ${\bf Q}_{{\mathcal T}_2 {\bf Y}}(\{\theta_1\}\times {\mathcal S}_1^2)$ is ``central'' in the collection of curves~$\{\{t_1\}\times{\mathcal S}_1^2 \vert\, t_1\in {\mathcal S}_1^1 \}$: call it a {\it poloidal equator}. Exchanging the roles of ${\mathcal S}_1^1$ and~${\mathcal S}_1^2$ similarly yields a {\it toroidal equator} ${\bf Q}_{{\mathcal T}_2 {\bf Y}}(\{{\mathcal S}_1^1\times \theta_2)$. Quantile regions then can take the form of a collection of toroidal caps with pole $\thetab$, of equatorial strips with  toroidal  or   poloidal equator ${\bf Q}_{{\mathcal T}_2 {\bf Y}}(\{{\mathcal S}_1^1\times \{\theta_2\})$ or ${\bf Q}_{{\mathcal T}_2 {\bf Y}}(\{\theta_1]\times{\mathcal S}_1^2)$. More precisely, these regions are obtained as follows.
\begin{enumerate}
\item[--] {\it (Toroidal caps)} In the flat torus representation of  ${\mathcal T}_2$ as an open square $(-\pi, \pi)\times (-\pi, \pi)$   centered at ${\bf F}_{{\mathcal T}_2 {\bf Y}}(\thetab)$, consider the squares ${\scriptstyle{\mathbb C}}(\tau)\coloneqq[-\tau^{1/2}\pi,\tau^{1/2}\pi]\times[-\tau^{1/2}\pi,\tau^{1/2}\pi]$;   the exponential of ${\scriptstyle{\mathbb C}}(\tau)$ at ${\bf F}_{{\mathcal T}_2 {\bf Y}}(\thetab)$ yields a region~${\scriptstyle{\mathbb C}}_{{\mathcal T}_2}(\tau)\coloneqq~\!\exp_{{\bf F}_{{\mathcal T}_2 {\bf Y}}(\thetab)}\big[ {\scriptstyle{\mathbb C}}(\tau)\big]$ of~${\mathcal T}_2$ with~${\rm U}_{_{{\mathcal T}_2}}$ probability content $\tau$ and a quantile region of order~$\tau$ as the {\it toroidal cap}~${{\mathbb C}}_{{\mathcal T}_2 {\bf Y}}(\tau)\coloneqq {\bf Q}_{{\mathcal T}_2 {\bf Y}}({\scriptstyle{\mathbb C}}_{{\mathcal T}_2}(\tau))$, which has ${\rm P}_{\bf Y}$ probability content $\tau$ and is homotopic to~$\thetab$.
\item[--] {\it (Toroidal equatorial strips)} Still in  the flat torus representation of  ${\mathcal T}_2$ centered at~${\bf F}_{{\mathcal T}_2 {\bf Y}}(\thetab)$, consider the bands   ${\scriptstyle{\mathbb C}}(\tau)\coloneqq[-\pi\tau ,\pi\tau]\times(-\pi,\pi)$. The image ~${\scriptstyle{\mathbb C}}_{{\mathcal T}_2}(\tau)\coloneqq~\exp_{{\bf F}_{{\mathcal T}_2 {\bf Y}}(\thetab)}\big[ {\scriptstyle{\mathbb C}}(\tau)\big]$ of ${\scriptstyle{\mathbb C}}(\tau)$    by the exponential  mapping at~${\bf F}_{{\mathcal T}_2 {\bf Y}}(\thetab)$ belongs to ${\mathcal T}_2$  and has~${\rm U}_{_{{\mathcal T}_2}}$ proba\-bility content $\tau$: define the quantile region of order $\tau$ of ${\rm P}_{\bf Y}$ as the {\it toroidal equatorial strip}~${{\mathbb C}}_{{\mathcal T}_2 {\bf Y}}(\tau)\coloneqq {\bf Q}_{{\mathcal T}_2 {\bf Y}}({\scriptstyle{\mathbb C}}_{{\mathcal T}_2}(\tau))$, which has ${\rm P}_{\bf Y}$ proba\-bility content $\tau$ and is homotopic to the toroidal equator  ${\bf Q}_{{\mathcal T}_2 {\bf Y}}(\{{\mathcal S}_1^1\times \{\theta_2\})$.   
\item[--] {\it (Poloidal equatorial strips)} Same as for the toroidal strips above, after interchanging the roles of the toroidal and poloidal equators.
\end{enumerate}
By construction, these three collections of quantile regions have ${\rm P}_{\bf Y}$ probability content $\tau$ irrespective of ${\rm P}_{\bf Y}$, as expected from quantile regions of order 
$\tau$. Along with these definitions of quantile regions and contours, the concepts of distribution and quantile functions on~${\mathcal T}_2$ developed by~\cite{HallinLiuRiemann24} thus satisfy all the desired properties (a)--(d); as  for (e),   the appropriate equivariance property is with respect to toroidal and poloidal symmetries and rotations.

\subsection{Conclusion}

The long quest for a satisfactory concept of multivariate quantile has generated many interesting and  useful developments, among which the various notions of depth, and the concepts of  $\rho$-, elliptical, spatial, and Oja quantiles. However, the measure-transportation-based center-outward quantile functions described in Section~\ref{Sec41} are the only ones to date that satisfy the two fundamental attributes of a quantile function:   full  charac\-terization of the underlying distribution ${\rm P}_{\bf Y}$ and  quantile regions of order $\tau$ with probability contents~$\tau$ irrespective of~${\rm P}_{\bf Y}$. They also nicely allow for meaningful nonparametric multiple-output and manifold-valued quantile regression and vector quantile autoregression \citep{Barrio25, Gonzetal26}. On the other hand, they are not naturally affine-equivariant, and their infinite-dimensional versions (there is no uniform distribution over the infinite-dimensional unit ball: see \cite{GHS26}) are problematic.

%\section{Multivariate distribution and quantile functions: empirical concepts}

%\subsection{}

\begin{funding}
The support of the Czech Science Foundation grant GA24-10078S and the Fonds Thelam of the Fondation Roi Baudouin is  gratefully acknowledged. 
\end{funding}
\bigskip

%\begin{thanks}
{\bf Acknowledgment. }
%Both authors acknowledge
The author would like to thank Dimitri Konen, Hang Liu, Gilles Mordant, Davy Paindaveine, and Mirek {\v Siman} for insightful discussions and comments.
%\end{thanks}

\bibliographystyle{imsart-nameyear}
\bibliography{MultQuant}

%\newpage

%\renewcommand{\thesection}{\Alph{section}}
%\setcounter{section}{0}

%\bibliographystyle{abbrvnat}
%\bibliography{MultQuant}
\end{document}